\documentclass[12pt]{article}
\usepackage{amsfonts}

\title{Fano's inequality is  false for a simple Cremona transformation
of  five-dimensional projective space }
\author{Marat Gizatullin\thanks{Department of Mathematics,
 Technical University Federico Santa Mar{\'{\i}}a,
Avenida Espa\~na, No.\ 1640, Casilla 110-V, Valpara{\'{\i}}so,
Chile. {\emph{e-mail}} \textsf{mgizatul@mat.utfsm.cl}}}
\date{14 of February 2002}

\newcommand{\PP}{\mathbb P}
\newcommand{\broken}{\dasharrow}
\newcommand{\mult}{\mathop{\rm mult}\nolimits}

\begin{document}
\maketitle

\par\noindent ABSTRACT. A Cremona transformation of
five-dimensional projective space is constructed. The degree of
the  transformation is 7.  The inequalities of Fano are not
fulfilled for this transformation.

\medskip
\maketitle \par\noindent Acknowledgement. I would like to thank
Miles Reid for helpful correspondence on the subject.
\section*{Introduction}
 It was shown in [3] and [4] that Fano's
inequalities are false for tree-dimensional projective space and
for three-dimensional quadric. Both the birational transformations
were of degree 13. Here we construct a simple five-dimensional
Cremona transformation   of degree 7 and not satisfying the
inequalities of Fano .
\medskip
\par  Let $X$ be a
non-singular $n$-dimensional variety such that
Pic$(X)=\mathbb{Z},$ and the anticanonical class $(-K_X) $ is
ample. If $(-K_X)=rH$ for a generator $H$ of the Picard group,
then $X$ is said to be Fano variety of index $r$ ( and of the
first kind ). According to the texts mentioned in survey [1],
 Fano's inequality is  the statement:
\begin{quote} \em For any birational
transformation,
 \[
 f\colon \ X \broken \ X
 \]
 defined by a linear system of  degree $d>1$
 (the degree is the
number defined by $f(mH)=dmH $ )
  there exists an irreducible subvariety $Y\subset X$ such that
 \[
 0\leq \dim Y \leq \dim X-2 ,
 \]
 \[
 \mult_Y(f(|H|))>{\frac{(\dim X-\dim Y -1)\cdot d }{r}} .
 \]
\end{quote}
\par
 Five-dimensional projective space $\PP^5$ has index 6, therefore
 one can reformulate Fano's inequality:
 \begin{quote} \em
For any Cremona transformation,
 \[
 f\colon \PP^5 \broken \PP^5
 \]
 defined by six homogeneous polynomials of the same degree $d>1$ and
without a common nonconstant factor,
 \[
 x_i'= f_i(x_0,x_1,x_2,x_3,x_4,x_5),\quad i=0,...,5,
 \]
  there exists an irreducible subvariety $Y\subset \PP^5$ such that
 \[
 0\leq \dim Y \leq 3 ,
 \]
 \[
 \mult_Y(f_i)> {\frac{(4-\dim Y )\cdot d }{6}} .
 \]
 for every $i=0,...,5$.
 \end{quote}
 \medskip
 \par
 The goal of my article is to show that these inequalities do not
take place for a certain Cremona transformation of degree 7. That
is, I write down the formulas for a Cremona transformation of
degree 7 such that for the forms $f_0,\dots,f_5$ defining the
transformation, for any surface $S \subset \PP^5,$   for any curve
$C\subset \PP^5$ and for any point $P\in\PP^5,$  one can see that
 \[
 {\min}_i(\mult_S(f_i))\le2 , \quad {\min}_i(\mult_C(f_i))\le3, \quad
{\min}_i(\mult_P(f_i))\le4.
 \]

\section*{Construction of the example}
Let us consider six homogeneous coordinates  for $\PP^5$ as the
normalized coefficients $x_{00},x_{01},x_{11},x_{02},x_{12},x_{22}
$ of a ternary quadratic form $F(T_0,T_1,T_2),$

 \[
 F(T_0,T_1,T_2)=x_{00}T_{0}^2+2x_{01}T_{0}T_{1}+x_{11}T_{1}^2
 \]
 \[
+2x_{02}T_{0}T_{2}+2x_{12}T_{1}T_{2}+x_{22}T_{2}^2.
\]
Let \[D=D(x_{00},x_{01},x_{11},x_{02},x_{12},x_{22})\] be the
discriminant  of the ternary quadric,
 \[
D={\left |\begin{array}{ccc}
  x_{00} & x_{01} & x_{02} \\
  x_{01} & x_{11} & x_{12} \\
  x_{02} & x_{12} & x_{22}   \end{array}\right | }.
 \]
 The set of double points of the cubic discriminant hypersurface
 consists of the points of the Veronese surface , these points
 correspond to the ternary quadrics which are perfect squares of
 linear forms.\par
 We fix  parameters $s,t $ and consider six following forms of degree 7.
\begin{eqnarray*}
 (f_{s,t})_{00} &=& x_{00}D^2, \\
 (f_{s,t})_{01} &=& x_{01}D^2, \\
 (f_{s,t})_{11}&=& x_{11}D^2,\\
 (f_{s,t})_{02} &=& x_{02}D^2+x_{01}x_{11}^3Ds+x_{00}^4Dt, \\
 (f_{s,t})_{12}&=&x_{12}D^2
 +x_{01}x_{00}^3Dt+x_{11}^4Ds,\\
 (f_{s,t})_{22}&=&x_{22}D^2+
2x_{12}x_{11}^3Ds+2x_{02}x_{00}^3Dt+x_{11}^7s^2+2x_{01}x_{11}^3x_{00}^3st
+x_{00}^7t^2.
 \end{eqnarray*}
These six forms define a two-parameter family of rational maps
 \[
 g_{s,t}\colon \PP^5\broken \PP^5.
 \]
 If $s=t=0$, all four forms have a common factor $D^2$. After cancelling
this, we see that $g_{0,0}$ is the identity transformation. For
our example, we need  nonzero values of $s,t$. If  one of the
parameters $s,t$ is not zero, then it is clear that the six forms
have no nonconstant common factor. Moreover,
 \[
 D\Big((f_{s,t})_{00},(f_{s,t})_{01},(f_{s,t})_{11},(f_{s,t})_{02},(f_{s,t})_{12},(f_{s,t})_{22}\Big)=D^{7}.
 \]
In fact, this identity expresses the invariance of the
discriminant under triangular transformation of variables
$T_0,T_1,T_2$. Using the latter identity, it is not hard to see
that
 \[
 (f_{-s,-t})_{ij}\Big((f_{s,t})_{00},(f_{s,t})_{01},(f_{s,t})_{11},(f_{s,t})_{02},(f_{s,t})_{12},(f_{s,t})_{22}\Big)=x_{ij}D^{16},
 \]
that is,
 \[
 g_{-s,-t}\circ g_{s,t}=\hbox{the identity transformation}.
 \]
 Thus $g_{s,t}$ is rationally invertible and is a Cremona transformation.
More generally,
 \[
 g_{s,t}\circ g_{p,q}=g_{s+p,t+q},
 \]
and we get a two-parameter group of Cremona transformations. These
transformations induce biregular automorphisms of an affine open
subset of the five-dimensional projective space, the complement of
the discriminant cubic hypersurface $D=0$. Indeed, the above
formula of the discriminant transformation proves this (moreover,
one can see below the exact calculation of the fundamental points
of such a transformation).
\bigskip
\par {\textbf{Remark.}}
The formulas for $g_{s,t}$ can be generalized for obtaining  an
infinite dimensional family of automorphisms of the complement to
the discriminant hypersurface.  The  construction of general
formulas resembles the trick used on page 8 of the Max-Planck
Institute preprint [2]. If $U(x,y,z)$ is a ternary cubic form,
$\phi_m(u,d),\psi_m(u,d)$ are binary forms of degree $m$,
\[ \phi=\phi_m(U(x_{00},x_{01},x_{11}),D),\]
\[ \psi=\psi_m(U(x_{00},x_{01},x_{11}),D),\]
then following transformation
\begin{eqnarray*}
 x_{00}^{\prime} &=& x_{00}D^{2m}, \\
 x_{01}^{\prime} &=& x_{01}D^{2m}, \\
 x_{11}^{\prime}&=& x_{11}D^{2m},\\
 x_{02}^{\prime} &=& x_{02}D^{2m}+x_{01}D^m\psi+x_{00}D^m\phi, \\
 x_{12}^{\prime}&=&x_{12}D^{2m}
 +x_{01}D^m\phi+x_{11}D^m\phi,\\
 x_{22}^{\prime}&=&x_{22}D^{2m}+
2x_{12}D^m\phi +2x_{02}D^m\psi+x_{11}\phi^2+2x_{01}\phi\psi
+x_{00}\psi^2.
 \end{eqnarray*}
is a Cremona transformation inducing a biregular automorphism of
the mentioned affine subset.
\bigskip
\par
Let us return to our two-parameter family. We fix  nonzero values
of the parameters $s,t$, for example, $s=t=1$, and consider the
corresponding Cremona transformation
\begin{eqnarray*}
 x_{00}^{\prime} &=& x_{00}D^2, \\
 x_{01}^{\prime} &=& x_{01}D^2, \\
 x_{11}^{\prime}&=& x_{11}D^2,\\
 x_{02}^{\prime} &=& x_{02}D^2+(x_{01}x_{11}^3+x_{00}^4)D, \\
 x_{12}^{\prime}&=&x_{12}D^2
 +(x_{01}x_{00}^3+x_{11}^4)D,\\
 x_{22}^{\prime}&=&x_{22}D^2+
2(x_{12}x_{11}^3+x_{02}x_{00}^3)D+(x_{11}^7+2x_{01}x_{11}^3x_{00}^3
+x_{00}^7).
 \end{eqnarray*}
 First of all, we find the points $P\in\PP^5$ where each form on the right
hand side has positive multiplicity (that is, the set of all
common zeros of these right hand sides, or, in other words, the
fundamental points of the transformation). The first right hand
side vanishes if \begin{itemize}\item either $x_{00}=0$,\item  or
$D=0$, \item or simultaneously $x_{00}=D=0$.\end{itemize}
\bigskip \par If $x_{00}=0$ but $D\ne 0$ then using the other five
formulas, one sees that for other five coordinates of a
fundamental point, the equalities
$x_{01}=x_{11}=x_{02}=x_{12}=x_{22}=0$ take place. This case is
not a point of $\PP^5$. \bigskip \par  The case $D=0$ but
$x_{00}\ne 0$ is  realizable.
 The points satisfying
 \[
 x_{00}\ne 0\quad  D=0,\quad  x_{11}^7+2x_{01}x_{11}^3x_{00}^3
+x_{00}^7=0
\] are fundamental, but they are of multiplicity 1
for one of three forms   $x_{02}',  x_{12}', x_{22}' $ at least.
Indeed, if such a point is double on the discriminant
hypersurface, then then it is non-singular on the  hypersurface
 of seventh degree.
If both the expressions  $(x_{01}x_{11}^3+x_{00}^4)$ ,\quad
$(x_{01}x_{00}^3+x_{11}^4)$ \quad  in
 $x_{02}',  x_{12}'$
vanish , and the point is fundamental , then the multiplicity of
$x_{22}'$ at the point is equal to 1.
\bigskip \par If  $x_{00}=D=0$ and the point is fundamental, then
also $x_{11}=0$ and either $x_{01}=0$ or
$(2x_{02}x_{12}-x_{01}x_{22})=0,$ or both the expressions vanish .
If the point is simple on the discriminant cubic and has
multiplicity more than 2 on every $x_{ij}', $ then all the
homogeneous coordinates vanish. If the point is double on $D=0$,
then $x_{01}=x_{02}=x_{12}=0,$ that is the only non-zero
coordinate of the point is $x_{22}. $ The multiplicity of $x_{22}'
$ at the point is 4.



\begin{thebibliography}{1000}
\bibitem[1]{1}
 Pukhlikov A.V. ,{\emph{Birational automorphisms of higher-dimensional algebraic varieties.}}
 Proceedings of
the International Congress of Mathematicians,  Berlin,  Vol. II,
(1998), 97-107.
\bibitem[2]{2}
 M.Gizatullin.
 {\emph{Examples of $m$-algebras}}, Max-Planck-Institut f{\"{u}}r
 Mathematik, Preprint Series, 2000 (50).
 \bibitem[3]{3}
 M.Gizatullin,
 {\emph{Fano's inequality is a mistake }}, E-preprint
math.arXiv.org , math AG/0202069 .
 \bibitem[4]{4}
 M.Gizatullin,
 {\emph{Fano's inequality is also false for
three-dimensional quadric}}, E-preprint math.arXiv.org , math
AG/0202117 .
\end{thebibliography}
  \end{document}